\documentclass[a4paper,10pt]{article}
\usepackage[top=3cm,bottom=3cm,left=2.5cm,right=2.5cm]{geometry}
\usepackage{amssymb}
\usepackage{graphicx}
\usepackage{amsmath,amsthm,amssymb,lineno}
\usepackage{latexsym}
\usepackage{graphicx,booktabs,multirow}
\usepackage{tikz}
\usetikzlibrary{decorations.pathreplacing}
\usetikzlibrary{intersections}
\usepackage{enumerate}
\usepackage{graphicx,booktabs,multirow}
\usepackage{appendix}
\newtheorem{theorem}{Theorem}
\newtheorem{lemma}[theorem]{Lemma}

\begin{document}


\title{The maximal tree with respect to the exponential of the second Zagreb index\thanks{This work is supported by the National Natural Science Foundation of China (11971164) and the Department of Education of Hunan Province (19A318).}}
\author{Mingyao Zeng, Hanyuan Deng\thanks{Corresponding author: hydeng@hunnu.edu.cn}\\
{\small College of Mathematics and Statistics, Hunan Normal University,}
 \\{\small  Changsha, Hunan 410081, P. R. China.}
}

\date{}
\maketitle

\begin{abstract}
The second Zagreb index is $M_2(G)=\sum_{uv\in E(G)}d_{G}(u)d_{G}(v)$. It was found to occur in certain approximate expressions of the total $\pi$-electron energy of alternant hydrocarbons and used by various researchers in their QSPR and QSAR studies. Recently the exponential of a vertex-degree-based topological index was introduced. It is known that among all trees with $n$ vertices, the exponential of the second Zagreb index $e^{M_2}$ attains its minimum value in the path $P_n$. In this paper, we show that $e^{M_2}$ attains its maximum value in the balanced double star with $n$ vertices and solve an open problem proposed by Cruz and Rada [R. Cruz, J. Rada, The path and the star as extremal values of vertex-degree-based topological indices among trees, MATCH Commun. Math. Comput. Chem. 82 (3) (2019) 715-732].

{\bf Keywords}: Exponential of the second Zagreb index; Maximal tree
\end{abstract}

\maketitle

\makeatletter
\renewcommand\@makefnmark%
{\mbox{\textsuperscript{\normalfont\@thefnmark)}}}
\makeatother

\baselineskip=0.25in

\section{Introduction}
\label{1}

The first Zagreb index $M_1$ and the second Zagreb index $M_2$ were introduced 50 years ago \cite{gt}. The Zagreb indices and their variants have been used to study molecular complexity, chirality, ZE-isomerism and heterosystems whilst the overall Zagreb indices exhibited a potential applicability for deriving multilinear regression models. Zagreb indices are also used by various researchers in their QSPR and QSAR studies. Mathematical properties of the Zagreb indices have also been studied. Readers can refer to the paper \cite{nkmt} and the cited literature. The second Zagreb index was found to occur in certain approximate expressions of the total $\pi$-electron energy of alternant hydrocarbons \cite{gt}. We encourage the reader to consult \cite{bkmr,dg,deng,dsa,eg,fs,jjt,ytl} for the historical background, computational techniques, and mathematical properties of the second Zagreb index.

For simple graph $G$ with edge set $E(G)$, the second Zagreb index of $G$ is defined as
$$M_2(G)=\sum_{uv\in E(G)}d_{G}(u)d_{G}(v)$$
where $d_{G}(v)$ is the degree of the vertex $v$ in $G$. It is a vertex-degree-based (VDB for short) topological index, also referred as bond incident degree index.

A formal definition of a VDB topological index is as follows. Let $\mathcal{G}_n$ be the set of graphs with $n$ non-isolated vertices. Consider the set
$$K =\{(i,j)\in \mathbf{N}\times \mathbf{N}: 1\leq i\leq j\leq n-1\}$$
and for a graph $G\in \mathcal{G}_n$, denote by $m_{i,j}(G)$ the number of edges in $G$ joining vertices of degree $i$ and $j$. A VDB topological index over $\mathcal{G}_n$ is a function $\varphi: \mathcal{G}_n\rightarrow \mathbf{R}$ induced by real numbers $\{\varphi_{i,j}\}_{(i,j)\in K}$ defined as
$$\varphi(G)=\sum_{(i,j)\in K}m_{i,j}(G)\varphi_{i,j}$$
for every $G\in \mathcal{G}_n$.

Many important topological indices are obtained from different choices of $\varphi_{i,j}$. For example, the first Zagreb index $M_1$ induced by numbers $\varphi_{i,j}=i+j$; the second Zagreb index $M_2$ induced by $\varphi_{i,j}=ij$; the Randi\'{c} index $\chi$ induced by $\varphi_{i,j}=\frac{1}{\sqrt{ij}}$, et al. For details on VDB topological indices, see \cite{cr,g,rada,tc}.

In order to study of the discrimination ability of topological indices,  Rada \cite{rada} introduced the exponential of a vertex-degree-based topological index. Given a vertex-degree-based topological index $\varphi$, the exponential of $\varphi$, denoted by $e^{\varphi}$, is defined as
$$e^{\varphi}(G)=\sum_{(i,j)\in K}m_{i,j}(G)e^{\varphi_{i,j}}$$
The extremal value problem of $e^{\varphi}$ over the set $\mathcal{T}_n$ of trees with $n$ vertices was initiated in \cite{cr}, and it was shown that $e^{M_1}$, $e^{M_2}$ attain their minimum value in the path $P_n$, $e^{\chi}$, $e^{H}$, $e^{GA}$, $e^{SC}$, $e^{AZ}$ attain their minimum value in the star $S_n$, $e^{M_1}$, $e^{ABC}$ attain their maximum value in the star $S_n$, $e^{H}$, $e^{GA}$, $e^{SC}$ attain their maximum value in the path $P_n$. In \cite{cmr}, it was shown that $e^{\chi}$ attains its maximum value in the path $P_n$. These results are summarized in Table \ref{tab-1}.

\begin{table}
\caption{Results on extremal trees for exponential of well known VDB topological indices.}\label{tab-1}
\begin{center}
\begin{tabular}{ccccccccc}
\hline
\quad& $e^{M_1}$ & $e^{M_2}$ & $e^{\chi}$ & $e^{H}$ & $e^{GA}$ &$e^{SC}$ & $e^{ABC}$ & $e^{AZ}$\\
\hline
min & $P_n$ & $P_n$ & $S_n$ & $S_n$ & $S_n$ & ${S_n}$ & $?$ & ${S_n}$\\
max & $S_n$ & $?$ & $P_n$ & $P_n$ & $P_n$ & ${P_n}$ & ${S_n}$ & $?$\\
\hline
\end{tabular}
\end{center}
\vskip -0.5cm
\end{table}

The maximum value of $e^{M_2}$ over $\mathcal{T}_n$ is still an open problem. In this paper, we prove that the maximum value of $e^{M_2}$ over $\mathcal{T}_n$ is attained in the balanced double star $S_{\lfloor\frac{n-2}{2}\rfloor, \lceil\frac{n-2}{2}\rceil}$, and solve an open problem proposed by Cruz and Rada \cite{cr}.

\section{Trees with maximum exponential second Zagreb index}

We first show in this section that in a maximal tree with respect to $e^{M_2}$ over $\mathcal{T}_n$, the distance between any pendant vertex and any vertex with the maximum degree in $T$ is at least 2

\begin{figure}[ht]
\begin{center}
  \includegraphics[width=10cm,height=3cm]{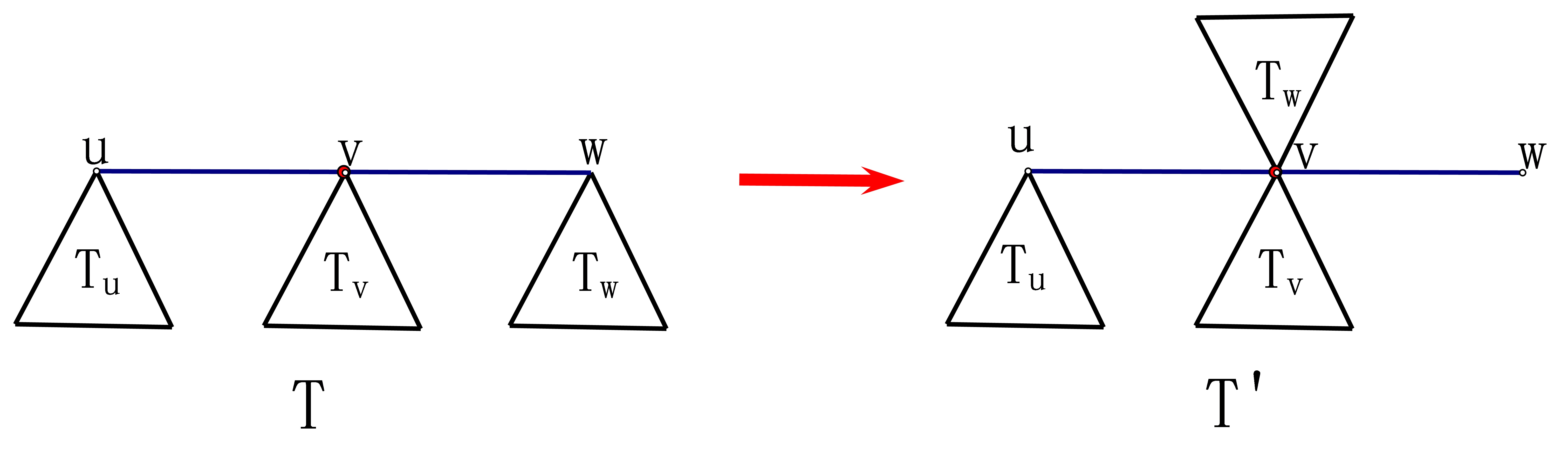}
 \end{center}
\vskip -0.5cm
\caption{Trees $T$ and $T'$ in Lemma \ref{l-2}}\label{fig-2}
\end{figure}

\begin{lemma}\label{l-2}
If $T$ is a maximal tree with respect to $e^{M_2}$ in $\mathcal{T}_n$, then the distance between a pendant vertex and a vertex with the maximum degree in $T$ is at most 2.
\end{lemma}
\textbf{Proof}. Otherwise, there is a vertex $u$ with the maximum degree in $T$ and a path $P=uvw$ such that $w$ is not a pendent vertex in $T$. Let $T_u$, $T_v$ and $T_w$ be the components of $T-uv-vw$ containing $u$, $v$ and $w$, respectively, see Figure \ref{fig-2}. Let $d(u)=\triangle$, $d(v)=s$, $d(w)=t$, where $\triangle\geq s\geq 2$, $\triangle\geq t\geq 2$. The set of neighbours of $v$ in $T$ is $N_T(v)=\{u,w,v_1,v_2,\cdots,v_{s-2}\}$ with $d_T(v_i)=x_i$ ($i=1,2,\cdots,s-2$), and the set of neighbours of $w$ is $N_T(w)=\{v,w_1,w_2,\cdots,w_{t-1}\}$ with $d(w_j)=y_j$ ($j=1,2,\cdots,t-1$). Let $T'=T-\{ww_1,ww_2,\cdots,ww_{t-1}\}+\{vw_1,vw_2,\cdots,vw_{t-1}\}$, then $d_{T'}(u)=\triangle$, $d_{T'}(v)=s+t-1$, $d_{T'}(w)=1$ and
\begin{align*}
e^{M_2}(T')-e^{M_2}(T)
&=e^{\triangle(s+t-1)}+\sum_{i=1}^{s-2}e^{x_i(s+t-1)}+\sum_{j=1}^{t-1}e^{y_j(s+t-1)}+e^{s+t-1}-(e^{\triangle s}+\sum_{i=1}^{s-2}e^{x_is}+\sum_{j=1}^{t-1}e^{y_jt}+e^{st})\\
&=(e^{\triangle(s+t-1)}-e^{\triangle s}-e^{st})+\sum_{i=1}^{s-2}(e^{x_i(s+t-1)}-e^{x_is})+\sum_{j=1}^{t-1}(e^{y_j(s+t-1)}-e^{y_jt})+e^{s+t-1}\\
&>0
\end{align*}

So, $e^{M_2}(T')>e^{M_2}(T)$ and $T$ is not a maximal tree with respect to $e^{M_2}$ in $\mathcal{T}_n$.
\hfill$\Box$

\begin{figure}[ht]
\begin{center}
  \includegraphics[width=5cm,height=5cm]{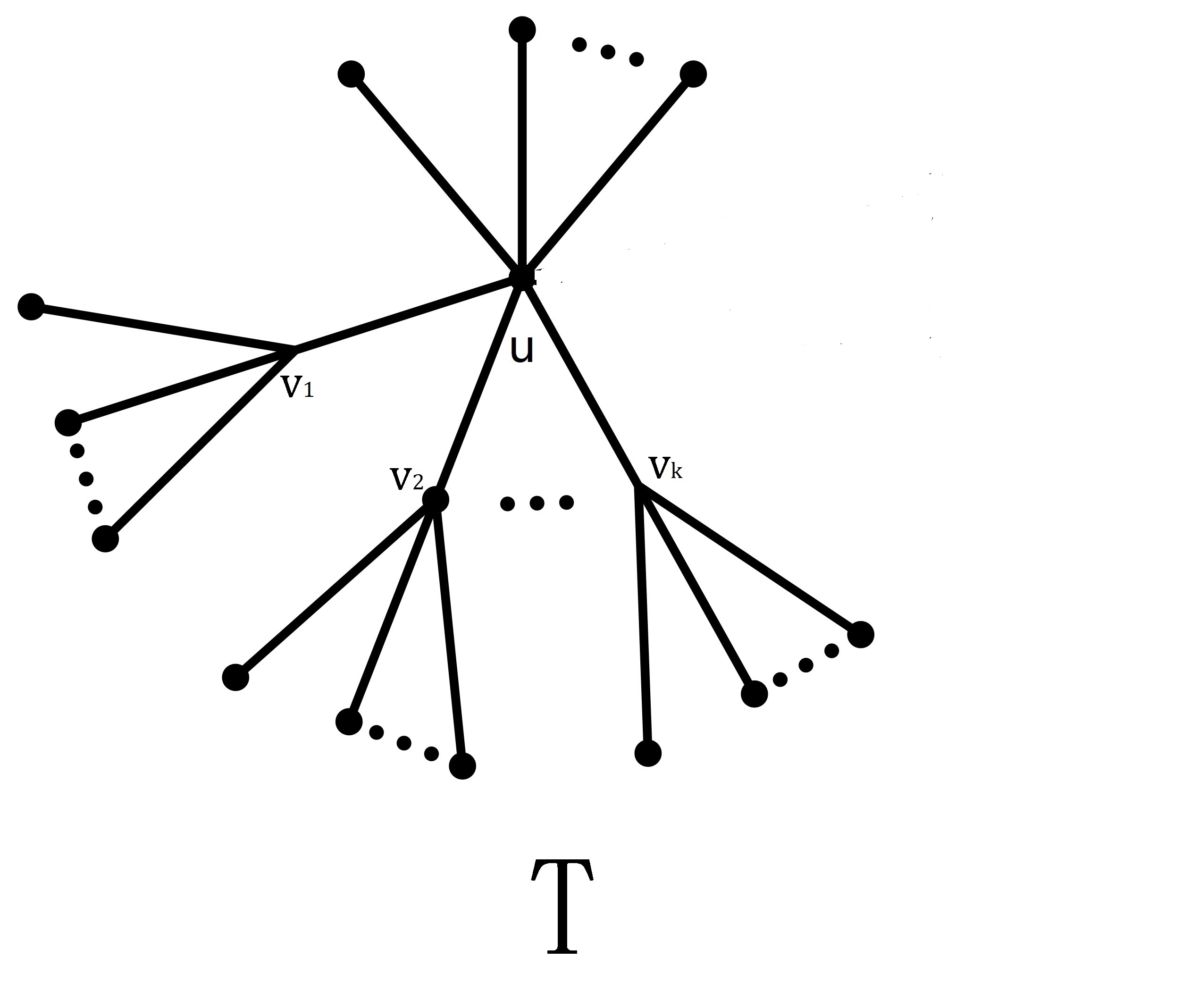}
 \end{center}
\vskip -0.5cm
\caption{A tree $T$}\label{fig-3}
\end{figure}

{\bf Remark}. Let $T\in \mathcal{T}_n$ ($n\geq 4$), in which all distances between any pendent vertex and any vertex with the maximum degree are at most 2, then $T$ has the form shown in Figure \ref{fig-3}, where $u$ is its unique vertex with the maximum degree if $k>1$, or $T$ is a double star if $k=1$.

\begin{figure}[ht]
\begin{center}
  \includegraphics[width=11cm,height=4.5cm]{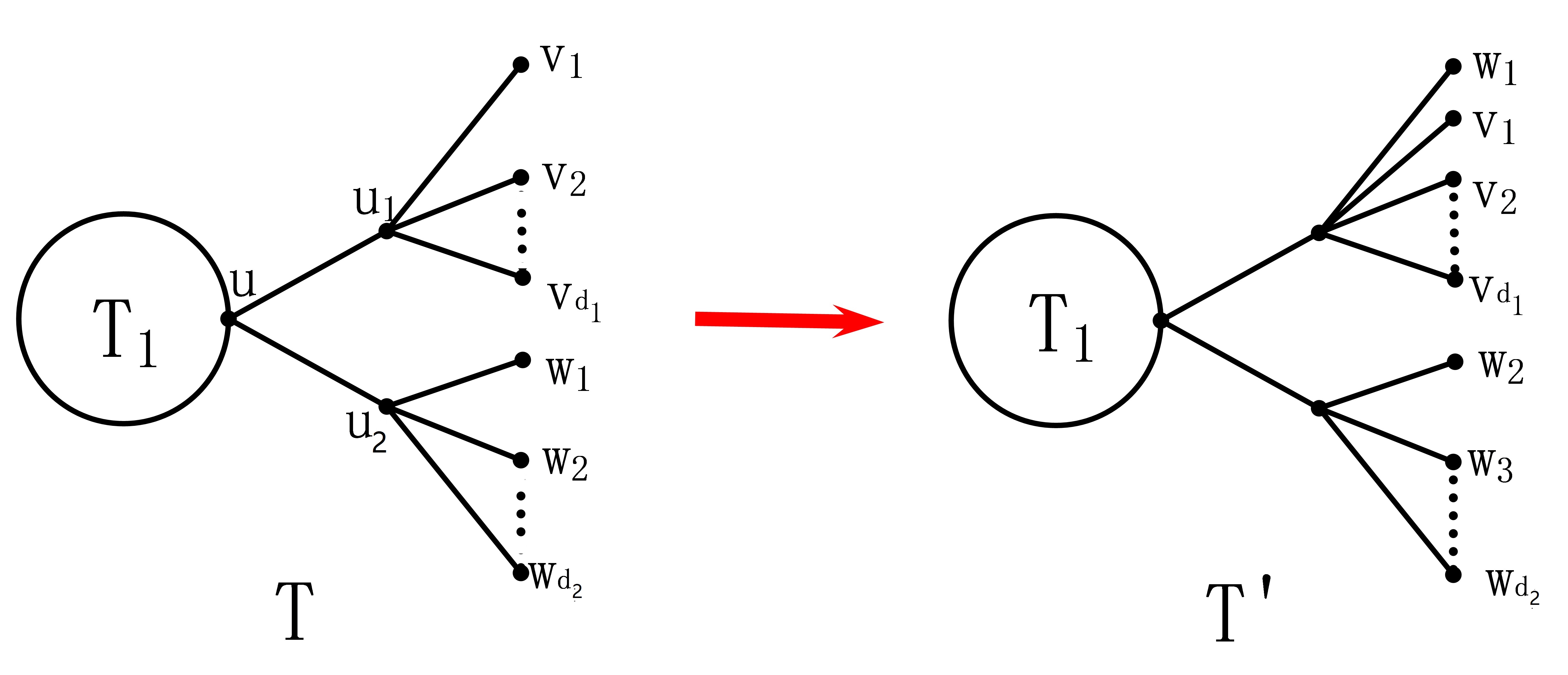}
 \end{center}
\vskip -0.5cm
\caption{The trees $T$ and $T'$ in Lemma \ref{l-3}.}\label{fig-4}
\end{figure}

\begin{lemma}\label{l-3}
Let $T$ and $T'$ be the trees on $n$ vertices given in Figure \ref{fig-4}, where $u$ is a vertex with the maximum degree and $T_1$ is a subtree of $T$. If $d_1\geq d_2\geq 1$, then $e^{M_2}(T)<e^{M_2}(T')$.
\end{lemma}
\textbf{Proof}.
Let $d_T(u)=\Delta$ be the maximum degree in $T$. The sets of neighbours of $u_1$ and $u_2$ in $T$ are $N_T(u_1)=\{u,v_1,\cdots,v_{d_1}\}$ and $N_T(u_2)=\{u,w_1,\cdots,w_{d_i}\}$, respectively, where $d(v_1)=d(v_2)=\cdots=d(v_{d_1})=1$, $d(w_1)=d(w_2)=\cdots=d(w_{d_2})=1$. Let $T'=T-\{u_2w_1\}+\{u_1w_1\}$, then $d_{T'}(u_1)=d_1+2$, $d_{T'}(u_2)=d_2$, and
\begin{align*}
e^{M_2}(T')-e^{M_2}(T)
=&e^{\Delta(d_1+2)}+(d_1+2)e^{(d_1+2)}+e^{\Delta d_2}+(d_2-1)e^{d_2}\\
&-[e^{\Delta(d_1+1)}+d_1e^{(d_1+1)}+e^{\Delta (d_2+1)}+d_2e^{(d_2+1)}]\\
=&[e^{\Delta(d_1+1)}\cdot e^{\Delta}-e^{\Delta(d_1+1)}]+[d_1e^{(d_1+1)}\cdot e-d_1e^{(d_1+1)}]\\
&+e^{(d_1+2)}+(e^{\Delta d_2}-e^{\Delta d_2}\cdot e^{\Delta})+(d_2e^{d_i}-d_2e^{d_2}\cdot e)-e^{d_2}\\
=&[e^{\Delta(d_1+1)}-e^{\Delta d_2}][e^{\Delta}-1]+[d_1e^{(d_1+1)}-d_2e^{d_2}](e-1)\\
&+[e^{(d_1+2)}-e^{{d_2}}]>0
\end{align*}
So, $e^{M_2}(T)<e^{M_2}(T')$.
\hfill$\Box$

In the following, we consider the double star $S_{x,y}$, it is a tree on $x+y+2$ vertices with exactly two non-pendent vertices of degrees $x+1$ and $y+1$, respectively.

\begin{lemma}\label{l-4}
Let $T=S_{x,y}$ be the double star on $n$ vertices, where $x+y=n-2$. If $T$ is a maximal tree with respect to $e^{M_2}$ among all double stars on $n$ vertices, then $|x-y|\leq 1$, i.e., $T$ is the balanced double star with $n$ vertices. And

$$e^{M_2}(S_{\lfloor\frac{n-2}{2}\rfloor, \lceil\frac{n-2}{2}\rceil})=\left\{
\begin{array}{ll}
  e^{\frac{n^2}{4}}+(n-2)e^{\frac{n}{2}}, & \mbox{if $n$ is even}; \\
  e^{\frac{n^2-1}{4}}+(\frac{n-3}{2})e^{\frac{n-1}{2}}+(\frac{n-1}{2})e^{\frac{n+1}{2}}, & \mbox{if $n$ is odd}.
\end{array}
\right.
$$
\end{lemma}
\textbf{Proof}.
Let $u$ and $v$ be two non-pendent vertices of the double star $T=S_{x,y}$, $d_T(u)=x+1$ and $d_T(v)=y+1$. Without loss of generality, we assume that $x\leq y$.
If $|x-y|>1$, let $T'=T-vv_1+uv_1$, where $v_1\neq u$ is a neighbour of $v$ in $T$, then $T'\simeq S_{x+1,y-1}$ and
\begin{align*}
e^{M_2}(T')-e^{M_2}(T)
&=[e^{(x+2)y}+(x+1)e^{(x+2)}+(y-1)e^{y}]-[e^{(x+1)(y+1)}+xe^{x+1}+ye^{y+1}]\\
&=e^{(x+2)y}-e^{(x+1)(y+1)}-ye^{y+1}+(x+1)e^{(x+2)}-xe^{x+1}+(y-1)e^{y}\\
&>e^{xy+2y}-e^{xy+x+y+1}-ye^{y+1}\\
&=e^{xy+y}[e^y-e^{x+1}]-ye^{y+1}\\
&\geq e^{xy+y}[e^{x+2}-e^{x+1}]-ye^{y+1}\\
&=e^{xy+y+x+1}[e-1]-ye^{y+1}\\
&>e^{xy+y+x+1}-ye^{y+1}\\
&=e^{y+1}[e^{x(y+1)}-y]>0
\end{align*}
So, $e^{M_2}(T')>e^{M_2}(T)$ and $T$ is not maximal with respect to $e^{M_2}$ among all double stars on $n$ vertices.
\hfill$\Box$

\begin{theorem}\label{t-1}
If $T\in \mathcal{T}_n$ and $T\not\simeq S_{\lfloor\frac{n-2}{2}\rfloor, \lceil\frac{n-2}{2}\rceil}$, $n\geq 4$, then $T$ is not maximal with respect to $e^{M_2}$ over $\mathcal{T}_n$.
\end{theorem}
\textbf{Proof}.
In fact, $e^{M_2}(S_{\lfloor\frac{n-2}{2}\rfloor, \lceil\frac{n-2}{2}\rceil})>e^{M_2}(S_n)$ for $n\geq 4$, the star $S_n$ on $n$ vertices is not maximal with respect to $e^{M_2}$ over $\mathcal{T}_n$.

If $T$ is not a double star, then by Lemma \ref{l-2}, we may assume that all distances between a pendant vertex and a vertex with the maximum degree in $T$ are at least 2 and $T$ has the form depicted in Figure \ref{fig-3}, where $k>1$. By Lemma \ref{l-3}, $T$ is not maximal.

So, $T$ is a double star and the result follows from Lemma \ref{l-4}.
\hfill$\Box$

\section*{References}

\bibliographystyle{elsarticle-num}
\bibliography{<your-bib-database>}

\end{document}